\documentclass[10pt]{article}
\usepackage{amssymb}
\usepackage{amsmath}
\usepackage{amsfonts}
\usepackage{hyperref}
\usepackage{graphicx}
\usepackage{amsmath,cite}
\usepackage{bookmark}
\usepackage { amssymb }
\usepackage{bm}
\usepackage{amsthm}

\usepackage[OT1]{fontenc}
\DeclareFontFamily{T1}{calligra}{}
\DeclareFontShape{T1}{calligra}{m}{n}{<->s*[1.44]callig15}{}
\DeclareMathAlphabet\mathcalligra   {T1}{calligra} {m} {n}
\DeclareMathAlphabet\mathzapf       {T1}{pzc} {mb} {it}
\DeclareMathAlphabet\mathchorus     {T1}{qzc} {m} {n}
\DeclareMathAlphabet\mathrsfso      {U}{rsfso}{m}{n}

\newcommand{\keywords}[1]{\par\smallskip\noindent\textbf{Keywords.} #1}
\newtheorem{theorem}{Theorem}

\newtheorem{corollary}{Corollary}

\newtheorem{example}{Example}

\newtheorem{lemma}{Lemma}

\begin{document}
\title{Generalized Fundamental Theorem of Calculus on a Regular Time Scale}
\author{Billur KAYMAK\c{C}ALAN$^{1}$ \quad Yeter ŞAHİNER$^{2}$\\[0.3em]
Şeyda YILDIZ$^{3}$ \quad Emre Deniz YILMAZ$^{4}$\\[1em]
\parbox{0.94\textwidth}{\centering\small
$^{1}$Faculty of Engineering, University of Turkish Aeronautical Association, Ankara, TÜRKİYE\\[0.2em]
$^{2}$Elementary Mathematics Teaching Programme, Hacettepe University, Ankara, TÜRKİYE\\[0.2em]
$^{3}$Faculty of Arts and Science, Department of Mathematics, TED University, Ankara, TÜRKİYE\\[0.2em]
$^{4}$Bilkent University, 06800 Bilkent, Ankara, TÜRKİYE}}
\date{}
\maketitle
\begin{abstract}
The main objectives of this paper can be summarized as, giving a fundamental theorem on a regular time scale as well as illustrating that applications are possible on a homogeneous time scale.
\end{abstract}

\section{Introduction}\label{sec1}

\noindent Time scale theory provides a framework that unifies discrete and continuous dynamical systems under a single structure \cite{hilger1}. Within this framework, it becomes possible to analyze both difference and differential equations in a unified manner.\\

\noindent A time scale (TS) is a closed subset of $\mathbb{R}$ that may include both discrete and continuous components \cite{Lak,bohner}. Typical examples include $\mathbb{R}$, $\mathbb{Z}$, or more structured sets such as grids and fractals. This approach broadens the classical understanding of time and makes it possible to investigate dynamical systems that cannot be completely characterized by either continuous or discrete formulations alone.\\

\noindent Overall, the concept of TS's is important because it provides a flexible and unified framework for studying dynamic systems with diverse temporal characteristics. It enables researchers to apply mathematical tools and techniques to analyze complex phenomena and develop deeper insights into the behavior of systems over time.\\

\noindent Among the various types of TS's studied in this theory, homogeneous time scales play a pivotal role due to their smooth and consistent structure, \cite{homogenous,control,homogenous1,sb} and references therein. A homogeneous time scale $\bm{\textbf{T}}$ is characterized by its uniformity and regularity, ensuring a constant graininess of time throughout its domain. This property makes homogeneous time scales particularly well-suited for mathematical modeling and analysis, as they provide a seamless integration of continuous and discrete dynamics and this means that the TS does not exhibit any abrupt changes in its structure or density. In a homogeneous time scale, the principles of calculus extend naturally, allowing for the development of differential and difference methods, integration techniques, and stability theory. The details concerning homogeneous as well as regular time scales are given Section 3.1 below.\\

\noindent The studies of homogeneous time scales \cite{homogenous,control} have led to significant advancements in various fields of mathematics and science. From dynamical systems and control theory \cite {control} to mathematical biology \cite{thomas} and economics \cite{atici1}, homogeneous time scales offer a structured approach to investigating systems in which continuous and discrete dynamics coexist. Moreover, the theory provides valuable insights into the fundamental properties of dynamic systems, shedding light on stability, controllability, and predictability.\\

\noindent The classical Fundamental Theorem of Calculus (FTC) on TS's \cite{Lak,bohner,bohner1,atici} generalizes the classical FTC. The classical FTC provides a direct link between integration and differentiation by showing that, for a continuous function $f$ on $[a,b]$, its definite integral can be computed using any anti-derivative $F$ as follows:
$$\displaystyle\int\limits_a^b  f(x)dx=F(b)-F(a).$$

\noindent In a similar manner, FTC can be adapted to functions defined on TS's which may include both continuous and discrete elements. Within this framework, the theorem has been formulated for both $\Delta$  \cite{Lak,bohner} and $\nabla$  integrals \cite{bohner1,atici,gurses}. These formulations provide a consistent relationship between differentiation and integration that reflects the structure of the underlying TS. However, such a result is not generally available for the $\Diamond_\alpha$ integral \cite{sheng1, sheng}.\\

\noindent Motivated by this limitation, the present study focuses on regular time scales and aims to establish a generalized version of FTC for the $\Diamond_\alpha$ integral. In addition, the structure of homogeneous time scales is examined, and it is shown how the obtained results can be applied in this setting. The results of this study offer new insights into diamond-alpha calculus and further elucidate its role in the context of TS theory.

\section{Preliminaries}

\noindent Typical examples of a TS $\mathbf{T}$ include the continuous case $\mathbf{T} = \mathbb{R} $, the discrete case 
$\mathbf{T} = \mathbb{Z}$, and more general structures such as $\mathbf{T}=h\mathbb{Z}$.\\

\noindent For each $s \in \mathbf{T}$, the forward and backward jump operators are defined by
\[
\sigma(s) = \inf \{ t \in \mathbf{T} \ \ |\ \ t > s \}, \quad 
\nu(s) = \sup \{ t \in \mathbf{T}  \ \ |\ \ t < s \}.
\]

\noindent The corresponding graininess functions are given by
\[
\mu(s) = \sigma(s) - s, \quad \eta(s) = s - \nu(s).
\]

\noindent A point $s \in \mathbf{T}$ is called right-scattered if $\sigma(s) > s$, and right-dense if $\sigma(s) = s$. 
Similarly, $s$ is called left-scattered if $\nu(s) < s$, and left-dense if $\nu(s) = s$.\\

\noindent A function $\mathbf{x} : \mathbf{T} \to \mathbb{R}$ is said to be delta differentiable at $s \in \mathbf{T}$ if there exists a number $\mathbf{x}^{\Delta}(s)$ satisfying the corresponding limit condition on the TS \cite{bohner1}. 
In particular, if $s$ is right-scattered, then
\begin{equation}\label{delta}
    \mathbf{x}^{\Delta}(s) = \frac{\mathbf{x}(\sigma(s)) - \mathbf{x}(s)}{\mu(s)},
\end{equation}

\noindent while for right-dense points it reduces to the classical derivative.\\

\noindent Similarly, the nabla derivative $\mathbf{x}^{\nabla}(s)$ is defined using the backward jump operator. If $s$ is left-scattered, then

\begin{equation}\label{nabla}
\mathbf{x}^{\nabla}(s) = \frac{\mathbf{x}(s) - \mathbf{x}(\nu(s))}{\eta(s)},
\end{equation}

\noindent and for left-dense points it coincides with the usual limit definition.\\

\noindent A function is called rd-continuous if it is continuous at right-dense points and admits finite left limits at left-dense points. The class of real rd-continuous functions defined on a TS $\mathbf{T}$ is denoted by $\mathbf{C}_{rd}(\mathbf{T};\mathbb{R})$. The notion of ld-continuity is defined analogously with $\mathbf{C}_{ld}(\mathbf{T};\mathbb{R})$. These function classes guarantee the existence of delta and nabla integrals, which are introduced via anti-derivatives similarly to the classical case.\\

\noindent For $0 \leq \alpha \leq 1$, the diamond-$\alpha$ derivative is defined, whenever both derivatives exist, by
\[
\mathbf{x}^{\Diamond_\alpha}(s) = \alpha \mathbf{x}^{\Delta}(s) + (1-\alpha)\mathbf{x}^{\nabla}(s).
\]

\noindent The corresponding diamond-$\alpha$ integral is defined as a convex combination of the delta and nabla integrals:
\[
\int_a^b \mathbf{x}(s)\,\Diamond_\alpha s 
= \alpha \int_a^b \mathbf{x}(s)\,\Delta s 
+ (1-\alpha)\int_a^b \mathbf{x}(s)\,\nabla s,
\]

\noindent for all $\mathbf{x}\in\mathbf{C}_{rl}(\mathbf{T};\mathbb{R})$, where $\mathbf{C}_{rl}(\mathbf{T};\mathbb{R}):=\mathbf{C}_{rd}(\mathbf{T};\mathbb{R})\cap \mathbf{C}_{ld}(\mathbf{T};\mathbb{R}).\medskip$

\noindent Unlike the classical setting, the relationship between the diamond-$\alpha$ derivative and integral is not straightforward, which motivates further investigation.\\

\noindent In general, $\Diamond_\alpha$ derivative of $\displaystyle\int\limits_a^t  \mathbf{x}(s)\Diamond_\alpha s$ with respect to $t$ is not equal to $\mathbf{x}(t)$. However, studies on the calculation of $\ \Bigg(\displaystyle\int\limits_a^t  \mathbf{x}(s)\Diamond_\alpha s\Bigg)^{\Diamond_\alpha}$ are gaining importance. It has been a matter of curiosity, for which TS's or functions defined on those TS's, the fundamental theorem is valid. With this motivation, we start our study by considering regular time scales.

\section{MAIN RESULTS}

\noindent In this study, some remarkable results are obtained for regular time scales (RTS). These results have been obtained by imposing restrictions on\\

\noindent \textbf{(a) } the functions defined on a RTS.\\

\noindent \textbf{(b) } the graininess functions of a RTS.\\

\noindent A TS $\bm{\textbf{T}}$ is referred to as \textbf{regular} if it is invariant under a double application of the forward jump operator and if the backward graininess evaluated at the forward jump restores the original point \cite{gurses}; that is,
\[
\sigma(\sigma(s)) = s \quad \text{and} \quad \nu(\sigma(s)) = s, \quad \forall s \in \bm{\textbf{T}}.
\]

\noindent It is easy to see that graininess functions are invertible on RTS's, that is, $\sigma^{-1}=\nu$ and $\nu^{-1}=\sigma$ and also, the following statements can be given on RTS's.\\

\begin{lemma} \label{lmm}
Let $\ \bm{\textbf{T}}$ be a RTS. Then the followings hold: 
\\

\noindent \textbf{(a) } $ \ \eta\big(\sigma(s)\big)=\mu(s)$ and $ \ \mu\big(\sigma(s)\big)=\eta(s)$.\\

\noindent \textbf{(b) } $ \ \sigma^\nabla(s)=\dfrac{\mu(s)}{\eta(s)}$ and $ \nu^\Delta(s)=\dfrac{\eta(s)}{\mu(s)}, $ and hence $\sigma^\nabla(s) \nu^\Delta(s)=1$.

\end{lemma}

\begin{proof}
\noindent \textbf{(a) }  $\mu(s)=\sigma(s)-s=\sigma(s)-\nu\big(\sigma(s)\big)=\eta\big(\sigma(s)\big) $ and the other is similar.\\

\noindent \textbf{(b) } Using  \eqref{nabla}, we obtain $\sigma^\nabla(s)=\dfrac{\sigma(s)-\sigma\big(\sigma(s)\big)}{\eta(s)}=\dfrac{\sigma(s)-s}{\eta(s)}=\dfrac{\mu(s)}{\eta(s)} $ and the other is similar.

\end{proof}

\begin{lemma}\label{lmm3} Let $\ \bm{\textbf{T}}$ be a RTS. Then the followings hold for a function $\mathbf{x}:\bm{\textbf{T}}\rightarrow \mathbb{R}$:\\

\noindent \textbf{(a) } $\Big[  \mathbf{x}(\nu(s))\Big]^\Delta=\nu^\Delta(s) \mathbf{x}^\Delta(\nu(s))\quad \quad $  \textbf{(b) } $\Big[  \mathbf{x}(\sigma(s))\Big]^\nabla=\sigma^\nabla(s) \mathbf{x}^\nabla(\sigma(s))$.\\

\end{lemma}
\begin{proof} \textbf{(a)}
    $$\Big[  \mathbf{x}(\nu(s))\Big]^\Delta=\dfrac{\mathbf{x}(\nu(\sigma(s)))-\mathbf{x}(\nu(s))}{\mu(s)}=\dfrac{\mathbf{x}(\sigma(\nu(s)))-\mathbf{x}(\nu(s))}{\eta(s)}\dfrac{\eta(s)}{\mu(s)}$$ then we get $$\Big[  \mathbf{x}(\nu(s))\Big]^\Delta=\nu^\Delta(s)\dfrac{\mathbf{x}(\sigma(\nu(s)))-\mathbf{x}(\nu(s))}{\mu(\nu(s))}=\nu^\Delta(s) \mathbf{x}^\Delta(\nu(s))$$ and the other is similar.\\
\end{proof}

\noindent Another interesting TS is the \textbf{homogeneous time scale}, in which $$\mu(s)=\eta(s)=c,\ c\in\mathbb{R}$$ for all $s$.\\

\noindent Every homogeneous time scale is regular. Therefore this means that there are some TS's which are regular, but not homogeneous. For example, the TS's $$\bm{\textbf{T}}=(-\infty,0] \bigcup \bigg\{\dfrac{1}{n} \quad\Big|\quad n\in\mathbb{N} \bigg\}\bigcup \bigg\{\dfrac{n}{n+1} \quad\Big|\quad n\in\mathbb{N} \bigg\} \bigcup\big[1,2\big]$$ and $\bm{\textbf{T}}=\overline{q^\mathbb{Z}}$ are both regular, but not homogeneous; however the TS's $\bm{\textbf{T}} = h\mathbb{Z}$ for $h > 0\medskip$ and $\bm{\textbf{T}} = \mathbb{R}$ are both homogeneous and regular.\\

\noindent Certain restrictions are necessary for functions defined on a TS to ensure the proper application of the Diamond-Alpha integral. First, the function must be regressive. Regressive functions guarantee that the forward jump operators remain within specific bounds on the TS. Additionally, the function should satisfy suitable smoothness conditions, which are particularly important for the existence and continuity of solutions in the context of dynamic equations. These constraints help preserve the mathematical consistency of the Diamond-Alpha integral and contribute to obtaining reliable results in practical applications. However these are necessery but not enough for $\Diamond_\alpha$ integral calculus. Therefore, the fundamental theorem is necessary in this context. In order to obtain the fundamental theorem, the need for functions defined on the TS or some restrictions that can be imposed on the TS's is becoming more essential. \\

\noindent We start our work with the restrictions that can be placed on functions defined on a regular time scale and give some examples. Next, we will investigate the calculation of the $\Diamond_\alpha$ integral of functions with some restrictions on the graininess functions of a regular time scale. \\

\subsection{Construction of Generalized Fundamental Theorem on a RTS} \label{ssec1}

\noindent Anti-derivative definitions have been given for delta and nabla integrals, but there is no such definition for diamond-alpha integrals in the literature. The following definition is given to fill this gap:\\

\noindent Let us define a family of functions $\mathcal{X}$ as a set $$ \mathcal{X}=\Big\{ \mathbf{x} \ \big| \ \mathbf{x}\in\mathbf{C}_{rl}(\bm{\textbf{T}};\mathbb{R}) \ \text{and} \ \mathbf{x}\big(\sigma(s)\big)+\mathbf{x}\big(\nu(s)\big)=2\mathbf{x}(s)\Big\}$$

\noindent $\bm{\textbf{T}}$ is a RTS; $a,s\in\bm{\textbf{T}}$ and $\bm{\mathsf{X}} : \bm{\textbf{T}} \rightarrow \mathbb{R}$. A function $\mathsf{X}\in\mathcal{X}$ is called a \textbf{$\Diamond_\alpha$ antiderivative of $\mathbf{x}\in\mathbf{C}_{rl}(\bm{\textbf{T}};\mathbb{R})$} provided $\mathsf{X}^{\Diamond_\alpha}(s)=\mathbf{x}(s)$ holds for $s\in\bm{\textbf{T}}$. We define the \textbf{$\Diamond_\alpha$ integral of $\mathbf{x}$} by  $$\displaystyle\int\limits_a^s \mathbf{x}(t)\Diamond_\alpha t=\bm{\mathsf{X}}(s)-\bm{\mathsf{X}}(a).$$

\begin{theorem}[\textbf{Generalized Fundamental Theorem on a RTS}]\label{thm2}
Suppose that $\bm{\textbf{T}}$ is a RTS and $\mathbf{x}\in \mathcal{X}$.

\noindent\textbf{Part I. } If $\mathsf{X}(s)=\displaystyle\int\limits_a^s \mathbf{x}(t)\Diamond_\alpha t$, then $\mathsf{X}^{\Diamond_\alpha}(s)=\mathbf{x}(s)$.\\

\noindent\textbf{Part II. } $\displaystyle\int\limits_a^b \mathbf{x}(s)\Diamond_\alpha s=\bm{\mathsf{X}}(b)-\bm{\mathsf{X}}(a)$.\\
\end{theorem}
\begin{proof}

\noindent\textbf{Part I.) }
Since $\bm{\textbf{T}}$ is regular, we know that $\sigma\big(\sigma(s)\big)=s=\nu\big(\sigma(s)\big)$. Using Lemma \ref{lmm} - \textbf{(a)}, we obtain

$$ \displaystyle\int\limits_{s}^{\sigma(s)} \mathbf{x}(t)\nabla t = \displaystyle\int\limits_{\nu\big(\sigma(s)\big)}^{\sigma(s)} \mathbf{x}(t)\nabla t = \eta\big(\sigma(s)\big)\mathbf{x}\big(\sigma(s)\big)=\mu(s) \mathbf{x}\big(\sigma(s)\big)$$ 
and 
$$ \displaystyle\int\limits_{\nu(s)}^{s} \mathbf{x}(t)\Delta t = \displaystyle\int\limits_{\nu(s)}^{\sigma\big(\sigma(s)\big)} \mathbf{x}(t)\Delta t = \mu\big(\sigma(s)\big)\mathbf{x}\big(\sigma(s)\big)=\eta(s)\mathbf{x}\big(\sigma(s)\big).\medskip\\$$
\begin{equation*}
\begin{array}{ll}
\mathsf{X}^{\Diamond_\alpha}(s) & =\alpha \bm{\mathsf{X}}^\Delta (s)+(1-\alpha)\bm{\mathsf{X}}^\nabla (s) \medskip\\
& =\alpha \Bigg\{\displaystyle\int\limits_a^s \mathbf{x}(t)\Diamond_\alpha t\Bigg\}^\Delta+(1-\alpha) \Bigg\{\displaystyle\int\limits_a^s \mathbf{x}(t)\Diamond_\alpha t\Bigg\}^\nabla \medskip\\
& =\dfrac{\alpha}{\mu(s)} \Bigg\{\displaystyle\int\limits_a^{\sigma(s)} \mathbf{x}(t)\Diamond_\alpha t-\displaystyle\int\limits_a^s \mathbf{x}(t)\Diamond_\alpha t\Bigg\} \medskip\\
&\quad+\dfrac{1-\alpha}{\eta(s)} \Bigg\{\displaystyle\int\limits_a^s \mathbf{x}(t)\Diamond_\alpha t-\displaystyle\int\limits_a^{\nu(s)} \mathbf{x}(t)\Diamond_\alpha t\Bigg\} \medskip\\
& =\dfrac{\alpha}{\mu(s)} \displaystyle\int\limits_s^{\sigma(s)} \mathbf{x}(t)\Diamond_\alpha t+\dfrac{1-\alpha}{\eta(s)}\displaystyle\int\limits_{\nu(s)}^s \mathbf{x}(t)\Diamond_\alpha t\medskip\\
& =\dfrac{\alpha}{\mu(s)} \Bigg[\alpha\displaystyle\int\limits_s^{\sigma(s)} \mathbf{x}(t)\Delta t+(1-\alpha) \displaystyle\int\limits_{s=\nu\big(\sigma(s)\big)}^{\sigma(s)} \mathbf{x}(t)\nabla t \Bigg]\medskip\\
&+\dfrac{1-\alpha}{\eta(s)}\Bigg[\alpha\displaystyle\int\limits_{\nu(s)}^{s=\sigma\big(\nu(s)\big)} \mathbf{x}(t)\Delta t+(1-\alpha) \displaystyle\int\limits_{\nu(s)}^s \mathbf{x}(t)\nabla t \Bigg] \medskip\\
& =\dfrac{\alpha}{\mu(s)} \Big[ \alpha\mu(s)\mathbf{x}(s)+(1-\alpha)\mu(s)\mathbf{x}\big(\sigma(s)\big) \Big]\medskip\\
&+\dfrac{1-\alpha}{\eta(s)} \Big[ \alpha\eta(s)\mathbf{x}\big(\nu(s)\big) +(1-\alpha)\eta(s)\mathbf{x}(s) \Big] \medskip\\
\end{array} 
\end{equation*}
then,
\begin{equation*}
\begin{array}{ll}
\mathsf{X}^{\Diamond_\alpha}(s) & =\alpha^2\mathbf{x}(s)+\Big(\alpha-\alpha^2\Big)\mathbf{x}\big(\sigma(s)\big) +\Big(\alpha-\alpha^2\Big)\mathbf{x}\big(\nu(s)\big) +(1-\alpha)^2\mathbf{x}(s) \medskip\\
& =\Big(2\alpha^2-2\alpha+1\Big)\mathbf{x}(s)+\Big(\alpha-\alpha^2\Big)\Big[\mathbf{x}\big(\sigma(s)\big) +\mathbf{x}\big(\nu(s)\big)\Big] \medskip\\
& =\Big(2\alpha^2-2\alpha+1\Big)\mathbf{x}(s)+2\Big(\alpha-\alpha^2\Big)\mathbf{x}(s) \medskip\\
& =\Big[(1-\alpha)^2+\alpha^2+2\Big(\alpha-\alpha^2\Big)\Big]\mathbf{x}(s) \medskip\\
& =\mathbf{x}(s). \medskip\\
\end{array} 
\end{equation*}

\noindent\textbf{Part II.) }    Let $G(s)=\displaystyle\int\limits_a^s \mathsf{X}^{\Diamond_\alpha}(t)\Diamond_\alpha t$ by Part-I. We have $G^{\Diamond_\alpha}(s)=\mathsf{X}^{\Diamond_\alpha}(s)$, so $G(s)-\mathsf{X}(s)=C$ for every $s\in\bm{\textbf{T}}$ where $C$ is an arbitrary constant.\\

   \noindent When $s=a$, we get $G(a)-\bm{\mathsf{X}}(a)=C$. Since $G(a)=0$, we have $C=-\bm{\mathsf{X}}(a)$. So, $G(s)-\mathsf{X}(s)=-\bm{\mathsf{X}}(a)$ is provided.\\

   \noindent When $s=b$, we get $G(b)-\bm{\mathsf{X}}(b)=-\bm{\mathsf{X}}(a)$, i.e. $G(b)=\bm{\mathsf{X}}(b)-\bm{\mathsf{X}}(a)$ or $$\displaystyle\int\limits_a^b \mathsf{X}^{\Diamond_\alpha}(s)\Diamond_\alpha s=\bm{\mathsf{X}}(b)-\bm{\mathsf{X}}(a).$$
 \end{proof}

\noindent Some examples of fundamental theorem can be given as follows:\\

\begin{example}
The TS $\bm{\textbf{T}}=\bigg\{ n^2 \ \Big| \ n>1,\ n\in\mathbb{N} \bigg\} $ is regular since  $\sigma\big(\sigma(s)\big)=\nu\big(\sigma(s)\big)=s$ in $s\in\bm{\textbf{T}}$. For this TS, $\sigma(s)=\big( \sqrt{s}-1\big)^2$, $\nu(s)=\big( \sqrt{s}+1\big)^2$, $\eta(s)=2\sqrt{s}+1$ and $\mu(s)=2\sqrt{s}-1$. The function $\bm{\mathsf{X}}(s)=\sqrt{s}$ satisfies  $\bm{\mathsf{X}}\big(\sigma(s)\big)+\bm{\mathsf{X}}\big(\sigma(s)\big)=2\bm{\mathsf{X}}(s)$. So, $\mathsf{X}\in\mathcal{X}$ and $\Diamond_\alpha$ antiderivative of $\mathbf{x}(s)=\dfrac{(1-2\alpha)+2\sqrt{s}}{4s-1}$, $\alpha\in(0,1)$ since $\mathsf{X}^{\Diamond_\alpha}(s)=\mathbf{x}(s)$. By definition, we see that 
$$\displaystyle\int\limits_a^s \dfrac{(1-2\alpha)+2\sqrt{t}}{4t-1}\ \Diamond_\alpha t=\sqrt{t}\ \Bigg|_a^s.$$  
\end{example} 
 
\begin{example}
The TS $\bm{\textbf{T}}=\bigg\{ \dfrac{1}{n} \ \Big| \ n\in\mathbb{N} \bigg\} $ is regular since  $\sigma\big(\sigma(s)\big)=\nu\big(\sigma(s)\big)=s$ with $\sigma(s)=\dfrac{s}{1-s}$, $\nu(s)=\dfrac{s}{1+s}$. The function $\bm{\mathsf{X}}(s)=\dfrac{1}{s}$ satisfies  $\bm{\mathsf{X}}\big(\sigma(s)\big)+\bm{\mathsf{X}}\big(\sigma(s)\big)=2\bm{\mathsf{X}}(s)$. So, $\mathsf{X}\in\mathcal{X}$ and $\Diamond_\alpha$ antiderivative of $\mathbf{x}(s)=\dfrac{(2\alpha-1)s-1}{s^2}$, $\alpha\in(0,1)$ since $\mathsf{X}^{\Diamond_\alpha}(s)=\mathbf{x}(s)$. By definition, we see that 
$$\displaystyle\int\limits_a^s \dfrac{(2\alpha-1)t-1}{t^2}\ \Diamond_\alpha t=\dfrac{1}{t}\ \Bigg|_a^s.$$  
\end{example}

\begin{example}
The TS $\bm{\textbf{T}}=\bigg\{ \dfrac{n}{n+1} \ \Big| \ n\in\mathbb{N} \bigg\} $ is regular since  $\sigma\big(\sigma(s)\big)=\nu\big(\sigma(s)\big)=s$ in $s\in\bm{\textbf{T}}$. On $\bm{\textbf{T}}$, $\sigma(s)=\dfrac{1}{2-s}$ and $\nu(s)=\dfrac{2s-1}{s}$. The function $\bm{\mathsf{X}}(s)=\dfrac{s+1}{s-1}$ satisfies  $\bm{\mathsf{X}}\big(\sigma(s)\big)+\bm{\mathsf{X}}\big(\sigma(s)\big)=2\bm{\mathsf{X}}(s)$. So, $\mathsf{X}\in\mathcal{X}$ and $\Diamond_\alpha$ antiderivative of $\mathbf{x}(s)=\dfrac{(4\alpha-2)s-4\alpha}{(s-1)^2}$, $\alpha\in(0,1)$ since $\mathsf{X}^{\Diamond_\alpha}(s)=\mathbf{x}(s)$. By definition, we have 
$$\displaystyle\int\limits_a^s \dfrac{(4\alpha-2)t-4\alpha}{(t-1)^2}\ \Diamond_\alpha t=\dfrac{t+1}{t-1}\ \Bigg|_a^s.$$   
\end{example}

\subsection{Calculation of $\Diamond_\alpha$ Integral}

\noindent In this part, the calculation of $\Diamond_\alpha$ integration is investigated over both regular and homogeneous time scales.\\

\begin{lemma} \label{lmm2}
For any RTS $\bm{\textbf{T}}$ satisfying $\mu(s)=k \eta(s)$, $k\in\mathbb{R}$, \\

\noindent\textbf{(a) } $ \ \displaystyle\int\limits_a^b \mathbf{x}^\Delta(s)\Delta s=\dfrac{1}{k}\Bigg.\mathbf{x}(\sigma(s))\bigg|_a^b\quad \quad $  \textbf{(b) } $\ \displaystyle\int\limits_a^b \mathbf{x}^\nabla(s)\nabla s=k\Bigg.\mathbf{x}(\nu(s))\bigg|_a^b \quad$ hold.
\end{lemma}

\begin{proof}
We give the proof of \textbf{(a)}. The proof of \textbf{(b)} is similar.\\

$\displaystyle\int\limits_a^b \mathbf{x}^\Delta(s)\nabla s=\displaystyle\int\limits_a^b\mathbf{x}^\nabla\big(\sigma(s)\big)\nabla s=\displaystyle\int\limits_a^b\mathbf{x}^\nabla\big(\sigma(s)\big) \sigma^\nabla(s)\dfrac{\eta(s)}{\mu(s)}\nabla s,
$ \\

\noindent since $\quad\sigma^\nabla(s)\dfrac{\eta(s)}{\mu(s)}=1$ by Lemma \ref{lmm}. Using Lemma \ref{lmm3}, we have\\

$\displaystyle\int\limits_a^b \mathbf{x}^\Delta(s)\nabla s=\dfrac{1}{k}\displaystyle\int\limits_a^b\bigg[\mathbf{x}\big(\sigma(s)\big)\bigg]^\nabla\nabla s=\dfrac{1}{k}\Bigg.\mathbf{x}(\sigma(s))\Bigg|_a^b.$

\end{proof}

\begin{theorem}\label{thm4}
If $\bm{\textbf{T}}$ is a RTS satisfying $\mu(s)=k\ \eta(s)$, $k\in\mathbb{R}$ and $\mathbf{x}: \bm{\textbf{T}} \rightarrow \mathbb{R}$ is a $\Diamond_\alpha$ differentiable function, then for any $a,b\in\bm{\textbf{T}}$

\begin{equation}\label{11}
\displaystyle\int\limits_a^b  \mathbf{x}^{\Diamond_\alpha}(s)\Diamond_\alpha s=\Bigg\{\Big(2\alpha^2-2\alpha+1\Big)\mathbf{x}(s)+\Big(\alpha-\alpha^2\Big)\Bigg[\dfrac{\mathbf{x}\big(\sigma(s)\big)}{k}+k\mathbf{x}\big(\nu(s)\big)\Bigg]\Bigg\}\Bigg|_a^b.\medskip\\    
\end{equation}

\end{theorem}

\begin{proof}

Using definitions of $\Diamond_\alpha$ derivative and $\Diamond_\alpha$ integral, we obtain

$ 
\begin{array}{ll}
\displaystyle\int\limits_a^b  \mathbf{x}^{\Diamond_\alpha}(s)\Diamond_\alpha s &=\displaystyle\int\limits_a^b\bigg[ \alpha \mathbf{x}^\Delta(s)+(1-\alpha) \mathbf{x}^\nabla(s)\bigg]\Diamond_\alpha s \\
& \\
& =\alpha\displaystyle\int\limits_a^b  \mathbf{x}^\Delta(s) \Diamond_\alpha s+(1-\alpha) \displaystyle\int\limits_a^b \mathbf{x}^\nabla(s)\Diamond_\alpha s\\
& \\
& =\alpha \left\{\alpha \displaystyle\int\limits_a^b \mathbf{x}^\Delta(s)\Delta s+(1-\alpha)\displaystyle\int\limits_a^b \mathbf{x}^\Delta(s)\nabla s\right\}\\
& \\
& +(1-\alpha)\left\{\alpha \displaystyle\int\limits_a^b \mathbf{x}^\nabla(s)\Delta s+(1-\alpha)\displaystyle\int\limits_a^b \mathbf{x}^\nabla(s)\nabla s\right\}\\
& \\
& =\alpha^2\displaystyle\int\limits_a^b \mathbf{x}^\Delta(s)\Delta s+(1-\alpha)^2 \displaystyle\int\limits_a^b \mathbf{x}^\nabla(s)\nabla s\\
& \\
& +\alpha(1-\alpha)\left\{ \displaystyle\int\limits_a^b \mathbf{x}^\Delta(s)\nabla s+ \displaystyle\int\limits_a^b \mathbf{x}^\nabla(s)\Delta s\right\}.\medskip\\
\end{array}
$ 

\noindent Using Lemma \ref{lmm2} for the last term, we obtain the equality \eqref{11}. 
\end{proof}

\noindent For homogeneous time scales, we rewrite Theorem \ref{thm4} as a corollary.

\begin{corollary}
    Let $\bm{\textbf{T}}$ be a homogeneous time scale, $\eta(s)=c=\mu(s)$, $s\in\bm{\textbf{T}}$. For any $\Diamond_\alpha$ differentiable function $\mathbf{x}: \bm{\textbf{T}} \rightarrow \mathbb{R}$, the following holds:
    \begin{equation}\label{2}
     \displaystyle\int\limits_a^b  \mathbf{x}^{\Diamond_\alpha}(s)\Diamond_\alpha s=\Bigg\{\Big(2\alpha^2-2\alpha+1\Big)\mathbf{x}(s)+\Big(\alpha-\alpha^2\Big) \Big[\mathbf{x}(s+c)+\mathbf{x}(s-c)\Big]\Bigg\}\Bigg|_a^b\medskip\\  
    \end{equation}
\end{corollary}

\noindent In the examples given in Subsection \ref{ssec1}, the TS's are regular but the condition of Theorem \ref{thm4} is not valid for each function. By using Theorem \ref{thm4}, the $\Diamond_\alpha$ integral of some functions can be calculated. Let us give some examples:

\begin{example}
$\bm{\textbf{T}}=\mathbb{Z}$ is both homogeneous and regular time scale  where $c=1$, $\mu(s)=1=\eta(s)$, $k=1$. Since the conditions of Theorem \ref{thm4} are satisfied for any function $\mathbf{x} : \mathbb{Z} \rightarrow \mathbb{R}$, we have
\begin{equation*}
  \displaystyle\int\limits_a^b  \mathbf{x}^{\Diamond_\alpha}(s)\Diamond_\alpha s=\Bigg\{\Big(2\alpha^2-2\alpha+1\Big)\mathbf{x}(s)+\Big(\alpha-\alpha^2\Big)\Big[\mathbf{x}(s+1)+\mathbf{x}(s-1)\Big]\Bigg\}\Bigg|_a^b.\medskip\\ 
\end{equation*}

\end{example}
\begin{example}
$\bm{\textbf{T}}=h\mathbb{Z}$ is both homogeneous and regular time scale  where $c=h$, $\mu(s)=h=\eta(s)$, $k=1$. Since the conditions of Theorem \ref{thm4} are satisfied for any function $\mathbf{x} : h\mathbb{Z} \rightarrow \mathbb{R}$, we have
\begin{equation*}
  \displaystyle\int\limits_a^b  \mathbf{x}^{\Diamond_\alpha}(s)\Diamond_\alpha s=\Bigg\{\Big(2\alpha^2-2\alpha+1\Big)\mathbf{x}(s)+\Big(\alpha-\alpha^2\Big)\Big[\mathbf{x}(s+h)+\mathbf{x}(s-h)\Big]\Bigg\}\Bigg|_a^b.\medskip\\ 
\end{equation*}
\end{example}
\begin{example}
$\bm{\textbf{T}}=\mathbf{q^{\mathbb{Z}^+}}\ (q>1)$ is not homogeneous but regular time scale. $\mu(s)=(q-1)s$ and $\eta(s)=\dfrac{(q-1)}{s}$. Since the condition of Theorem \ref{thm4} is satisfied for $k=\dfrac{1}{q}$, then for any function $\mathbf{x} : \mathbf{\mathbf{q^{\mathbb{Z}^+}}} \rightarrow \mathbb{R}$, we have

\begin{equation*}
  \displaystyle\int\limits_a^b  \mathbf{x}^{\Diamond_\alpha}(s)\Diamond_\alpha s=\Bigg\{\Big[\alpha^2+(1-\alpha)^2\Big]\mathbf{x}(s)+\alpha(1-\alpha)\Bigg[q\ \mathbf{x}\big(qs\big)+\dfrac{1}{q}\mathbf{x}\bigg(\dfrac{s}{q}\bigg)\Bigg]\Bigg\}\Bigg|_a^b.\medskip\\ 
\end{equation*} 
\end{example}

\section{Integration of Derivative of a Function with $\Diamond_\alpha$}

\begin{lemma} 
Suppose that $\mathbf{T}$ is a RTS and the functions $\mathbf{x}$, $\mathbf{x}^\Delta$ and $\mathbf{x}^\nabla $ are continuous on $\mathbf{T}^\kappa_\kappa$. Then, $$ \eta(s)\Big(\mathbf{x}^\Delta(s)\Big)^\nabla=\mu(s)\Big(\mathbf{x}^\nabla(s)\Big)^\Delta. $$
\end{lemma}
\begin{proof}
Since $\mathbf{T}$ is regular, we have 
\begin{equation}\label{a1}
    \mathbf{x}^\nabla(s)=\mathbf{x}^\Delta\big(\nu(s)\big),\qquad \mathbf{x}^\Delta(s)=\mathbf{x}^\nabla\big(\sigma(s)\big)\medskip\\
\end{equation}  and 
\begin{equation}\label{a2}
   \mathbf{x}\big(\nu(s)\big)=\mathbf{x}(s)-\eta(s)\mathbf{x}^\nabla(s),\qquad \mathbf{x}\big(\sigma(s)\big)=\mathbf{x}(s)+\mu(s)\mathbf{x}^\Delta(s).\medskip\\
\end{equation}
After taking $\mathbf{x}^\Delta(s)$ and $ \mathbf{x}^\nabla(s) $ in place of the function $\mathbf{x}(s)$ in \eqref{a2}, respectively, and using \eqref{a1}, we get
\begin{equation*}\label{1c}
    \mathbf{x}^\nabla(s)=\mathbf{x}^\Delta\big(\nu(s)\big)=\mathbf{x}^\Delta(s)-\eta(s)\Big(\mathbf{x}^\Delta(s)\Big)^\nabla\medskip\\
\end{equation*}
and
\begin{equation*}\label{1d}
    \mathbf{x}^\Delta(s)=\mathbf{x}^\nabla\big(\sigma(s)\big)=\mathbf{x}^\nabla(s)+\mu(s)\Big(\mathbf{x}^\nabla(s)\Big)^\Delta\medskip\\
\end{equation*}
Next, $$ \mathbf{x}^\Delta(s)-\mathbf{x}^\nabla(s)=\mu(s)\Big(\mathbf{x}^\Delta(s)\Big)^\nabla=\eta(s)\Big(\mathbf{x}^\nabla(s)\Big)^\Delta$$ So, $$\eta(s)\Big(\mathbf{x}^\Delta(s)\Big)^\nabla=\mu(s)\Big(\mathbf{x}^\nabla(s)\Big)^\Delta$$
\end{proof}

\begin{corollary}\label{b1}
 If a TS $\mathbf{T}$ is homogeneous and the functions $\mathbf{x}$, $\mathbf{x}^\Delta$ and $\mathbf{x}^\nabla$ are continuous on $\mathbf{T}^\kappa_\kappa$, then $ \Big(\mathbf{x}^\nabla(s)\Big)^\Delta=\Big(\mathbf{x}^\Delta(s)\Big)^\nabla$, $\ \sigma^\nabla(s)=1$ and $\ \nu^\Delta(s)=1$.
\end{corollary}

\begin{lemma}\label{lmm3} If $\mathbf{T}$ is a RTS and $\mathbf{x} : \mathbf{T} \rightarrow \mathbb{R}$ is $\Delta-$ or $\nabla-$ differentiable function, then we have
    $\Big[  \mathbf{x}(\nu(s))\Big]^\Delta=\mathbf{x}^\Delta(\nu(s))\nu^\Delta(s)$ and     $\Big[  \mathbf{x}(\sigma(s))\Big]^\nabla=\mathbf{x}^\nabla(\sigma(s))\sigma^\nabla(s)\medskip$ hold.
\end{lemma}
\begin{proof}
\begin{align*}
\Big[\mathbf{x}(\nu(s))\Big]^\Delta
&=\dfrac{\mathbf{x}(\nu(\sigma(s)))-\mathbf{x}(\nu(s))}{\mu(s)}\\
&=\dfrac{\mathbf{x}(\sigma(\nu(s)))-\mathbf{x}(\nu(s))}{\eta(s)}\dfrac{\eta(s)}{\mu(s)}\\
&=\dfrac{\mathbf{x}(\sigma(\nu(s)))-\mathbf{x}(\nu(s))}{\mu(\nu(s))}\nu^\Delta(s)\\
&=\mathbf{x}^\Delta(\nu(s))\nu^\Delta(s)
\end{align*}
and the other is similar.
\end{proof}

\begin{lemma} \label{lmm2}
If $\mathbf{T}$ is a homogeneous time scale and $\mathbf{x} : \mathbf{T} \rightarrow \mathbb{R}$ is $\Delta-$ or $\nabla-$ differentiable function, we get\\

\noindent\textbf{(i) } $ \ \displaystyle\int\limits_a^b \mathbf{x}^\Delta(s)\nabla s=\Bigg.\mathbf{x}(\sigma(s))\bigg|_a^b\quad \quad $  \textbf{(ii) } $\ \displaystyle\int\limits_a^b \mathbf{x}^\nabla(s)\Delta s=\Bigg.\mathbf{x}(\nu(s))\bigg|_a^b \quad$ hold.
\end{lemma}

\begin{proof}
We give the proof of \textbf{(i)}. The proof of \textbf{(ii)} is similar.\\

$\displaystyle\int\limits_a^b \mathbf{x}^\Delta(s)\nabla s=\displaystyle\int\limits_a^b\mathbf{x}^\nabla\big(\sigma(s)\big)\nabla s=\displaystyle\int\limits_a^b\mathbf{x}^\nabla\big(\sigma(s)\big) \sigma^\nabla(s)\nabla s,
$ \\

\noindent since $\ \sigma^\nabla(s)=1$. Using Lemma \ref{lmm3}, we have\\

$\displaystyle\int\limits_a^b \mathbf{x}^\Delta(s)\nabla s=\displaystyle\int\limits_a^b\Bigg[\mathbf{x}\big(\sigma(s)\big)\Bigg]^\nabla\nabla s=\Bigg.\mathbf{x}(\sigma(s))\Bigg|_a^b$

\end{proof}

\begin{theorem}\label{thm2}
    If $\mathbf{T}$ is a homogeneous time scale and $\mathbf{x} : \mathbf{T} \rightarrow \mathbb{R} $ is a $\Diamond_\alpha$-differentiable function, then we obtain
    \begin{equation*}\label{2}
      \Bigg(\displaystyle\int\limits_a^t  \mathbf{x}(s)\Diamond_\alpha s\Bigg)^{\Diamond_\alpha} =\displaystyle\int\limits_a^t  \mathbf{x}^{\Diamond_\alpha}(s)\Diamond_\alpha s\medskip\\  
    \end{equation*}

\end{theorem}

\begin{proof}
By [an exploration on.. 2006-Sheng], it is known that
\begin{equation}\label{2.1}
     \Bigg(\displaystyle\int\limits_a^t  \mathbf{x}(s)\Diamond_\alpha s\Bigg)^{\Diamond_\alpha} =\Bigg\{\Big(2\alpha^2-2\alpha+1\Big)\mathbf{x}(s)+\alpha(1-\alpha)\Bigg[\mathbf{x}\big(\sigma(s)\big)+\mathbf{x}\big(\nu(s)\big)\Bigg]\Bigg\}\Bigg|_a^t\medskip\\  
    \end{equation}
\begin{equation*}
\begin{aligned}
\int_a^t \mathbf{x}^{\Diamond_\alpha}(s)\Diamond_\alpha s
&=\int_a^t\bigl[\alpha\mathbf{x}^\Delta(s)+(1-\alpha)\mathbf{x}^\nabla(s)\bigr]\Diamond_\alpha s\\[0.5em]
&=\alpha\int_a^t\mathbf{x}^\Delta(s)\Diamond_\alpha s
+(1-\alpha)\int_a^t\mathbf{x}^\nabla(s)\Diamond_\alpha s\\[0.5em]
&=\alpha\biggl\{\alpha\int_a^t\mathbf{x}^\Delta(s)\Delta s
+(1-\alpha)\int_a^t\mathbf{x}^\Delta(s)\nabla s\biggr\}\\
&\quad+(1-\alpha)\biggl\{\alpha\int_a^t\mathbf{x}^\nabla(s)\Delta s
+(1-\alpha)\int_a^t\mathbf{x}^\nabla(s)\nabla s\biggr\}\\[0.5em]
&=\alpha^2\int_a^t\mathbf{x}^\Delta(s)\Delta s
+(1-\alpha)^2\int_a^t\mathbf{x}^\nabla(s)\nabla s\\
&\quad+\alpha(1-\alpha)\biggl\{\int_a^t\mathbf{x}^\Delta(s)\nabla s
+\int_a^t\mathbf{x}^\nabla(s)\Delta s\biggr\}
\end{aligned}
\end{equation*}

\noindent According to Lemma \ref{lmm2}, 
\begin{equation}\label{2.2}
  \displaystyle\int\limits_a^t  \mathbf{x}^{\Diamond_\alpha}(s)\Diamond_\alpha s=\Bigg\{\Big(2\alpha^2-2\alpha+1\Big)\mathbf{x}(s)+\alpha(1-\alpha)\Bigg[\mathbf{x}\big(\sigma(s)\big)+\mathbf{x}\big(\nu(s)\big)\Bigg]\Bigg\}\Bigg|_a^t\medskip\\  
\end{equation}

\noindent Considering the equalities \eqref{2.1} and \eqref{2.2} together, we obtain 

\begin{equation*}
    \Bigg(\displaystyle\int\limits_a^t  \mathbf{x}(s)\Diamond_\alpha s\Bigg)^{\Diamond_\alpha} =\displaystyle\int\limits_a^t  \mathbf{x}^{\Diamond_\alpha}(s)\Diamond_\alpha s
\end{equation*}
\end{proof}

\noindent \textbf{Author Contribution Statements} The authors equally contributed to the paper.\\

\noindent \textbf{Declaration of Competing Interests} The authors declare that they have no known competing financial interests or personal relationships that could have appeared to influence the work reported in this paper.\\

\label{son}


\begin{thebibliography}{1} 

 \bibitem{atici1} Atici, F.M., Biles, D.C., Lebedinsky, A.T., An application of time scales to economics, {\it Math Comput. Model} (2006), 43(7–8) 718–726

 \bibitem{atici} Atici, F. M., Guseinov, G. Sh., On Green's functions and positive solutions for boundary value problems on time scales. {\it J. Comput. Appl. Math.} \textbf{141}, (2002), 75-99.

 \bibitem{homogenous} Bartosiewicz, Z., Kotta, Ü.,Pawłuszewicz, E., Wyrwas, M., Algebraic formalism of differential one-forms for nonlinear control systems on time scales, {\it Proc. Estonian Acad. Sci. Phys. Math.}, (2007), 56, 3, 264–282.

 \bibitem{control} Bartosiewicz, Z., Pawłuszewicz, E., Kotta, Ü., Wyrwas, M.,  Control systems on regular time scales and their differential rings, {\it Math. Control Signals Syst.} 22 (2011), 185-201.

 \bibitem{homogenous1} Bartosiewicz, Z., Belikov, J., Kotta, Ü., Wyrwas, M., State feedback linearization of nonlinear control systems on homogeneous time scales, {\it Nonlinear Analysis: Hybrid Systems}, 31 (2019), 69–85.

 \bibitem{sb} Bayeğ, S., Mert, R., Generalized Hukuhara diamond-alpha derivative of fuzzy valued functions on time scales, {\it Commun. Fac. Sci. Univ. Ank. Ser. A1 Math. Stat.}, Vol 74, no. 1, pp. 103–116.

 \bibitem{bohner1} Bohner, M., Peterson, A., Advances in Dynamic Equations on Time Scales, (2003), Birkhäuser, Boston.

 \bibitem{bohner} Bohner, M., Peterson, A., Dynamic Equations on Time Scales. An Introduction With Applications, (2001), Birkh\"{a}user Boston, Inc., Boston, MA.
 
\bibitem{gurses} G\"{u}rses, M., Guseinov, G. Sh., Silindir B., Integrable equations on time scales, {\it J. Math. Phys.} 46(11) (2005), 113510, 22.

\bibitem{hilger1} Hilger, S., Analysis on measure chains—a unified approach to continuous and discrete calculus, Results in Mathematics, 18 (1990), 18–56.

\bibitem{zey} Kayar Z., Kaymak\c{c}alan B. and Pelen N. N., Diamond-alpha Bennett-Leindler type dynamic inequalities, {\it Mathematical Methods in the Applied Sciences}, Vol 45, 5 (2022), 2797–2819; doi:10.1002/mma.7955.

\bibitem{Lak} Lakshmikantham, V., Sivasundaram, S., Kaymakçalan, B., Dynamic systems on measure chains, Vol 370 (1996). Kluwer Academic Publishers, MIA, Netherlands.
  
\bibitem{mozy} Mozyrska, D., Torres, Delfim F. M.. A study of diamond-alpha dynamic equations on regular time scales.  {\it Afr. Diaspora J. Math. (N.S.).} (2009), 8(1), 35-47.

\bibitem{sheng} Rogers, Jr. J. W., Sheng, Q., Notes on the diamond-$\alpha$ dynamic derivative on time scales, {\it J. Math. Anal. Appl. }, (2007), 326(1), 228-241.

\bibitem{sheng1} Sheng, Q., Fadag, M., Henderson, J., Davis, J. M., An exploration of combined dynamic derivatives on time scales and their applications, {\it  Nonlinear Anal. Real World Appl.} 7(3) (2006), 395-413.

\bibitem{thomas} Thomas, D., Weedermann, M., Billings, L., Hoffacker, J., When to spray: a time scale calculus approach to controlling the impact of west nile virus. {\it Ecol Soc} Washington-Allen RA (2009) 14(2).

\end{thebibliography}
\end{document}